# Selection from a stable box

ALEXANDER AUE[1], ISTVÁN BERKES[2] and LAJOS HORVÁTH[3]

[1]*Department of Mathematical Sciences, Clemson University, O-324 Martin Hall, Clemson, SC 29634-0975, USA. E-mail:* alexaue@clemson.edu
[2]*Department of Statistics, Graz University of Technology, Steyrergasse 17/IV, A-8010 Graz, Austria. E-mail:* berkes@tugraz.at
[3]*Department of Mathematics, University of Utah, 155 South 1440 East, Salt Lake City, UT 84112-0090, USA. E-mail:* horvath@math.utah.edu

Let $\{X_j\}$ be independent, identically distributed random variables. It is well known that the functional CUSUM statistic and its randomly permuted version both converge weakly to a Brownian bridge if second moments exist. Surprisingly, an infinite-variance counterpart does not hold true. In the present paper, we let $\{X_j\}$ be in the domain of attraction of a strictly $\alpha$-stable law, $\alpha \in (0, 2)$. While the functional CUSUM statistics itself converges to an $\alpha$-stable bridge and so does the permuted version, provided both the $\{X_j\}$ and the permutation are random, the situation turns out to be more delicate if a realization of the $\{X_j\}$ is fixed and randomness is restricted to the permutation. Here, the conditional distribution function of the permuted CUSUM statistics converges in probability to a random and nondegenerate limit.

## 1. Introduction

Let $X, X_1, \ldots, X_n, \ldots$ be independent, identically distributed random variables. CUSUM-based procedures, frequently used to test for time homogeneity of an underlying random phenomenon, are functionals of the process

$$Z_n(t) = \frac{1}{s_n \sqrt{n}} \sum_{j=1}^{\lfloor nt \rfloor} (X_j - \bar{X}_n), \qquad t \in [0, 1], \tag{1.1}$$

where $\lfloor \cdot \rfloor$ denotes integer part and

$$\bar{X}_n = \frac{1}{n} \sum_{j=1}^{n} X_j \quad \text{and} \quad s_n^2 = \frac{1}{n-1} \sum_{j=1}^{n} (X_j - \bar{X}_n)^2$$

are the sample mean and sample variance, respectively. The asymptotic properties of $\{Z_n(t)\}$ are immediate consequences of Donsker's invariance principle (see Billingsley (1968)). Denote by $\overset{\mathcal{D}[0,1]}{\longrightarrow}$ convergence in the space $\mathcal{D}[0,1]$ of cadlag functions equipped with the Skorokhod $J_1$-topology. We state the following fundamental theorem.

**Theorem 1.1.** *If* $EX^2 < \infty$, *then*

$$Z_n(t) \overset{\mathcal{D}[0,1]}{\longrightarrow} B(t) \qquad (n \to \infty),$$

*where* $\{B(t) : t \in [0,1]\}$ *is a Brownian bridge.*

One of the main problems in testing for a homogeneous data generating process is that of obtaining approximations for the critical values. It is well known, however, that the convergence of functionals of $Z_n(t)$ to their limit distributions is often quite slow (sometimes even the convergence speed is unknown), raising the question of whether Theorem 1.1 is reasonably applicable in the case of small sample sizes. To avoid these complications, along with bootstrap methods, the permutation principle has been suggested to simulate critical values. For a recent, comprehensive survey of this topic, we refer to Hušková (2004). Let $\pi = (\pi(1), \ldots, \pi(n))$ be a random permutation of $(1, \ldots, n)$ which is independent of $\mathbf{X} = (X_1, X_2, \ldots)$. The permuted version of $Z_n(t)$ is then defined by

$$Z_{n,\pi}(t) = \frac{1}{s_n\sqrt{n}} \sum_{j=1}^{\lfloor nt \rfloor} (X_{\pi(j)} - \bar{X}_n), \qquad t \in [0,1]. \tag{1.2}$$

One can think of $Z_{n,\pi}(t)$ as a CUSUM process whose summands are chosen from $X_1, \ldots, X_n$ by drawing from a box without replacement. The following result establishes the limiting behavior of $Z_{n,\pi}(t)$ and hence shows that the permutation method works to simulate critical values. Also, note that empirical evidence suggests that the convergence to the limit is faster than for $Z_n(t)$.

**Theorem 1.2.** *If* $EX^2 < \infty$, *then, for almost all realizations of* $\mathbf{X}$,

$$Z_{n,\pi}(t) \overset{\mathcal{D}[0,1]}{\longrightarrow} B(t) \qquad (n \to \infty),$$

*where* $\{B(t) : t \in [0,1]\}$ *is a Brownian bridge.*

Theorem 1.2 states that, for almost all realizations of the data, the permutation method provides the asymptotically correct critical values. A proof can be found in Billingsley (1968), pages 209–214.

In this paper we are interested in the asymptotics of the analogues of $Z_n(t)$ and $Z_{n,\pi}(t)$ in the infinite-variance case $EX^2 = \infty$. To this end, let $\alpha \in (0,2]$ and recall that a random variable $\xi_\alpha$ is strictly $\alpha$-stable if its characteristic function has the form

$$\phi_\alpha(t) = \begin{cases} \exp(-ct^2/2), & \text{if } \alpha = 2, \\ \exp\left(-c|t|^\alpha \left[1 - \mathrm{i}\beta\,\mathrm{sgn}(t)\tan\left(\dfrac{\pi}{2}\alpha\right)\right]\right), & \text{if } \alpha \in (0,1) \cup (1,2), \ \beta \in [-1,1], \\ \exp(-c|t|), & \text{if } \alpha = 1, \end{cases}$$

for some $c > 0$. This definition includes the normal law ($\alpha = 2$) and the Cauchy law ($\alpha = 1$), where the latter is moreover assumed to be symmetric. Note that $\mathrm{E}[\xi_\alpha]$ exists only for $\alpha \in (1, 2]$ and $\mathrm{E}[\xi_\alpha^2]$ only for $\alpha = 2$. Since scaling the random variables $X_n$ does not change the normalized partial sum processes in (1.5) and (1.6) below, for our purposes it will suffice to consider the case $c = 1$. The class of stable distributions has become more and more important in theory as well as in applications. Stable laws appear in a natural way in areas such as radio engineering, electronics, biology and economics (see Zolotarev (1986), Chapter 1). For other expositions on $\alpha$-stable random variables and processes, we refer to Samorodnitsky and Taqqu (1994) and Bingham *et al.* (1987).

In what follows, let $\alpha \in (0, 2)$ and consider independent, identically distributed random variables $\{X_j\}$ which are in the domain of attraction of a strictly $\alpha$-stable law. In the case $\alpha \neq 1$, this is equivalent to saying that there is a function $L(x)$, slowly varying at infinity, such that

$$\frac{1}{n^{1/\alpha}L(n)} \sum_{j=1}^{n} (X_j - \mu_\alpha) \xrightarrow{\mathcal{D}} \xi_\alpha \qquad (n \to \infty), \tag{1.3}$$

where $\xi_\alpha$ is a strictly $\alpha$-stable random variable and $\mu_\alpha = 0$ if $\alpha \in (0, 1)$. We say that $L$ is slowly varying at infinity if $\lim_{t \to \infty} L(ct)/L(t) = 1$ for all $c > 0$. In the case $\alpha = 1$, we will assume (1.3) with $\mu_\alpha = 0$, excluding nonlinear centering sequences which would prevent the weak convergence of the process $A_n(t)$ in (1.5) below. If (1.3) holds, then there is a $p \in [0, 1]$ such that the tail probabilities of $X$ and $|X|$ satisfy the relations

$$\lim_{y \to \infty} \frac{P\{X > y\}}{P\{|X| > y\}} = p \quad \text{and} \quad \lim_{y \to \infty} \frac{P\{X < -y\}}{P\{|X| > y\}} = q, \tag{1.4}$$

where $q = 1 - p$ and $\beta = 2p - 1$; see Feller (1966) and Petrov (1995) for details.

To define the counterparts of $Z_n(t)$ and $Z_{n,\pi}(t)$ in the infinite-variance case, introduce the process

$$A_n(t) = \frac{1}{T_n} \sum_{j=1}^{\lfloor nt \rfloor} (X_j - \bar{X}_n), \qquad t \in [0, 1], \tag{1.5}$$

and its permuted version

$$A_{n,\pi}(t) = \frac{1}{T_n} \sum_{j=1}^{\lfloor nt \rfloor} (X_{\pi(j)} - \bar{X}_n), \qquad t \in [0, 1], \tag{1.6}$$

where $T_n = \max_{j \le n} |X_j|$ and $\pi = (\pi(1), \ldots, \pi(n))$ is a random permutation of $(1, \ldots, n)$, independent of $\mathbf{X}$. In contrast to the finite-variance case, $T_n$ is a natural norming factor in (1.5) and (1.6) since $\sum_{j=1}^{n}(X_j - \mu_\alpha)/T_n$ has a nondegenerate limit distribution.

One might now expect $\{A_n(t)\}$ and $\{A_{n,\pi}(t)\}$ to follow an asymptotic pattern similar to that of $\{Z_n(t)\}$ and $\{Z_{n,\pi}(t)\}$. Consequently, given a realization of $\mathbf{X}$, one might want to use the permutation principle to simulate critical values that have a higher accuracy in the case of small sample sizes. Surprisingly, this turns out to be impossible. While a version of Theorem 1.1 holds true for $\{A_n(t)\}$, an adaptation of Theorem 1.2 cannot be given for almost all realizations of $\mathbf{X}$. On the contrary, we will show that, conditioned on $\mathbf{X}$, the distribution function of the permuted process $\{A_{n,\pi}(t)\}$ converges to a random and nondegenerate limit for every fixed $t$. On the other hand, there is an averaging effect: if both $\mathbf{X}$ and $\pi$ are assumed to be random, then $A_{n,\pi}(t)$ and $A_n(t)$ converge weakly to the same limit, namely to an $\alpha$-stable bridge.

The paper is organized as follows. In Section 2, we will precisely formulate our results. Their proofs will be given in Section 3.

## 2. Results

First, we study the limiting behavior for the sequence $\{A_n(t) : t \in [0,1]\}$. Its weak convergence can be derived from the joint convergence of

$$\left\{ \frac{1}{n^{1/\alpha}L(n)} \sum_{j=1}^{\lfloor nt \rfloor} (X_j - \mu_\alpha) : t \in [0,1] \right\} \quad \text{and} \quad \frac{T_n}{n^{1/\alpha}L(n)}.$$

To deal with these two sequences, we are going to apply the method developed in LePage *et al.* (1981), where it was proven that any stable distribution can be represented as an infinite sum of random variables constructed from partial sums of exponential random variables. More specifically, let $\{E_j\}$ be a sequence of independent, identically distributed exponential random variables having expected value 1 and define the sequence $\{Z_j\}$ by

$$Z_j = (E_1 + \cdots + E_j)^{-1/\alpha}, \qquad j \ge 1. \tag{2.1}$$

Additionally, let $\{\delta_j\}$ be independent, identically distributed random variables, independent of $\{E_j\}$, satisfying $P\{\delta_j = 1\} = p$ and $P\{\delta_j = -1\} = q$, with $p$ and $q$ specified in (1.4). Then, define

$$\eta = \begin{cases} \displaystyle\sum_{k=1}^{\infty} \delta_k Z_k, & \alpha \in (0,1], \\ \displaystyle\sum_{k=1}^{\infty} (\delta_k Z_k - (p-q)\mathrm{E}[Z_k I\{0 < Z_k \le 1\}]), & \alpha \in (1,2). \end{cases} \tag{2.2}$$

LePage *et al.* (1981), Theorem 1, showed that the sums in (2.2) converge with probability one and

$$\eta \overset{\mathcal{D}}{=} \xi_\alpha / c_1 \qquad \text{with some constant } c_1. \tag{2.3}$$

Consider a random vector $(W_\alpha(t), \mathcal{Z})$ such that $W_\alpha(t)$ is a strictly $\alpha$-stable process, $\mathcal{Z} \overset{\mathcal{D}}{=} Z_1$ and the joint distribution of $(W_\alpha(t), \mathcal{Z})$ is defined by the requirement that for any $0 = t_1 < \cdots < t_K = 1$, we have

$$(W_\alpha(t_2) - W_\alpha(t_1), W_\alpha(t_3) - W_\alpha(t_2), \ldots, W_\alpha(t_K) - W_\alpha(t_{K-1}), \mathcal{Z})$$
$$\overset{\mathcal{D}}{=} \left( (t_2 - t_1)^{1/\alpha} \eta_1, (t_3 - t_2)^{1/\alpha} \eta_2, \ldots, (t_K - t_{K-1})^{1/\alpha} \eta_{K-1}, \max_{1 \le j \le K-1} (t_{j+1} - t_j)^{1/\alpha} \mathcal{Z}_j \right),$$

where $(\eta_1, \mathcal{Z}_1), \ldots, (\eta_{K-1}, \mathcal{Z}_{K-1})$ are independent, identically distributed random vectors, distributed as $(\eta, \mathcal{Z})$. The following theorem is essentially known; see Kasahara and Watanabe (1986), Section 9.

**Theorem 2.1.** *If (1.3) holds, then*

$$\left( \frac{1}{n^{1/\alpha} L_1(n)} \sum_{j=1}^{\lfloor nt \rfloor} (X_j - \mu_\alpha), \frac{T_n}{n^{1/\alpha} L_1(n)} \right) \overset{\mathcal{D}^2[0,1]}{\longrightarrow} (W_\alpha(t), \mathcal{Z}) \qquad (n \to \infty),$$

*where $L_1(x) = c_1 L(x)$ with $c_1$ given in (2.3). The distribution of the vector $(W_\alpha(t), \mathcal{Z})$ is defined above.*

The asymptotic behavior of $\{A_n(t)\}$ is now an immediate consequence of Theorem 2.1 and the continuous mapping theorem.

**Corollary 2.1.** *If (1.3) holds, then*

$$A_n(t) \overset{\mathcal{D}[0,1]}{\longrightarrow} \frac{1}{\mathcal{Z}} B_\alpha(t) \qquad (n \to \infty),$$

*where, for $t \in [0,1]$, $B_\alpha(t) = W_\alpha(t) - t W_\alpha(1)$.*

The process $\{B_\alpha(t) : t \in [0,1]\}$ is sometimes called an $\alpha$-stable bridge. Very little is known about the distributions of functionals of $B_\alpha(t)/\mathcal{Z}$ and therefore it is hard to determine critical values needed to construct asymptotic test procedures. It would hence be beneficial to apply the permutation method to obtain critical values for functionals of $A_n(t)$. However, as the subsequent series of theorems shows, Theorem 1.2 does not have an infinite-variance counterpart.

Consider the process $\{A_{n,\pi}(t)\}$ defined in (1.6) and let both $\mathbf{X}$ and $\pi$ be random. We then obtain the following result.

**Theorem 2.2.** *If (1.3) holds, then*

$$A_{n,\pi}(t) \xrightarrow{\mathcal{D}[0,1]} \frac{1}{\mathcal{Z}} B_\alpha(t) \qquad (n \to \infty), \tag{2.4}$$

*where, for $t \in [0,1]$, $B_\alpha(t) = W_\alpha(t) - tW_\alpha(1)$.*

Obviously, Theorem 2.2 is much weaker than Theorem 1.2 since it is only valid for the average of the realizations. However, a stronger result, holding true for almost all realizations of **X**, cannot be proved. Let $P_\mathbf{X}$ and $E_\mathbf{X}$ denote conditional probability and expected value given **X**, respectively.

In order to state the main result of our paper, we introduce some further notation. Let $t \in (0,1)$ and, on some probability space, consider sequences $\{S_j : j \geq 1\}$, $\{S_j^* : j \geq 1\}$, $\{\delta_j(t) : j \geq 1\}$, $\{\delta_j^*(t) : j \geq 1\}$ such that these sequences are mutually independent, $\{S_j : j \geq 1\}$, $\{S_j^* : j \geq 1\}$ are partial sums of independent exponential random variables with parameter 1 and $\{\delta_j(t) : j \geq 1\}$, $\{\delta_j^*(t) : j \geq 1\}$ are both sequences of independent random variables taking the values 1 and 0 with probabilities $t$ and $1-t$, respectively. Let **S** be the $\sigma$-algebra generated by $\{S_j, S_j^* : j \geq 1\}$.

**Theorem 2.3.** *Let $t \in (0,1)$. If (1.3) holds, then*

$$P_\mathbf{X}\{A_{n,\pi}(t) \leq x\} \xrightarrow{\mathcal{D}} P_\mathbf{S}\{R(t) \leq x\} \qquad (n \to \infty)$$

*for any real $x$, where*

$$R(t) = \frac{1}{M}\left[-w_1 \sum_{j=1}^\infty \frac{1}{S_j^{1/\alpha}}(\delta_j(t) - t) + w_2 \sum_{j=1}^\infty \frac{1}{(S_j^*)^{1/\alpha}}(\delta_j^*(t) - t)\right]$$

*with $w_1 = q^{1/\alpha}$, $w_2 = p^{1/\alpha}$, where $p$ and $q$ are defined by (1.4) and*

$$M = \max\left\{\frac{w_1}{S_1^{1/\alpha}}, \frac{w_2}{(S_1^*)^{1/\alpha}}\right\}. \tag{2.5}$$

Theorem 2.3 immediately implies that $P_\mathbf{X}\{A_{n,\pi}(t) \leq x\}$ cannot converge to $P\{B_\alpha(t)/\mathcal{Z} \leq x\}$ with probability one. Note that Theorem 2.3 is the permutation analogue of Athreya (1987), which showed that the conditional distribution of the appropriately normalized bootstrap sample mean converges in distribution to a nondegenerate random variable. The limits in Athreya (1987) and Theorem 2.3 are different due to the sampling with replacement in the bootstrap case and without replacement in the permutation case.

Theorem 2.3 describes the limiting conditional distribution of $A_{n,\pi}(t)$ for any fixed $t$. Its proof can be modified to yield the asymptotic conditional distribution of vectors $(A_{n,\pi}(t_1), \ldots, A_{n,\pi}(t_r))$ for any choices of $0 < t_1 < \cdots < t_r < 1$. To do so, we must extend the definition of the sequences $\{\delta_j(t)\}$ and $\{\delta_j^*(t)\}$. Let $\{U_j\}$ and $\{U_j^*\}$ be independent, identically distributed random variables, uniform on $[0, 1]$, which are independent of $\{S_j\}$ and $\{S_j^*\}$, and let

$$\delta_j(t) = I\{U_j \leq t\} \quad \text{and} \quad \delta_j^*(t) = I\{U_j^* \leq t\} \qquad (j \geq 1). \tag{2.6}$$

Therein, $I\{A\}$ is the indicator function of a set $A$. Note that for a single $t$, this definition coincides with the definition in terms of Bernoulli variables given above. The latter, however, do not carry any information on the joint behavior for a collection $t_1, \ldots, t_r$. We arrive at the following theorem.

**Theorem 2.4.** *Let* $0 < t_1 < \cdots < t_r < 1$. *If (1.3) holds, then, as* $n \to \infty$,

$$P_{\mathbf{X}}\{A_{n,\pi}(t_1) \leq x_1, \ldots, A_{n,\pi}(t_r) \leq x_r\} \xrightarrow{\mathcal{D}} P_{\mathbf{S}}\{R(t_1) \leq x_1, \ldots, R(t_r) \leq x_r\}$$

*for any real* $x_1, \ldots, x_r$, *where* $R(t)$ *is defined in Theorem 2.3 with the* $\delta_j(t), \delta_j^*(t)$ *from (2.6).*

The proof of Theorem 2.3 will show that $A_{n,\pi}(t)$ depends only on the very large and very small order statistics. These, however, are asymptotically independent, explaining why the infinite sums defining $R(t)$ in Theorems 2.3 and 2.4 are independent. In addition, $\delta_j(t)$ is the limit of an indicator variable which determines whether or not the $j$th smallest order statistic is among the first $\lfloor nt \rfloor$ permuted observations. The same reasoning applies to $\delta_j^*(t)$ and the location of the $j$th largest order statistic in the permutation $\pi$.

The previous results concern the asymptotic properties of $\{A_n(t)\}$ and $\{A_{n,\pi}(t)\}$, but it is possible to consider modifications of these processes which involve replacing the normalization $T_n$ by a more general sequence

$$T_n^{(\nu)} = \left( \sum_{j=1}^n |X_j - \bar{X}_n|^\nu \right)^{1/\nu} \qquad \text{with some } \nu > \alpha.$$

The corresponding CUSUM processes are then defined by

$$A_n^{(\nu)}(t) = \frac{1}{T_n^{(\nu)}} \sum_{j=1}^{\lfloor nt \rfloor} (X_j - \bar{X}_n), \qquad A_{n,\pi}^{(\nu)}(t) = \frac{1}{T_n^{(\nu)}} \sum_{j=1}^{\lfloor nt \rfloor} (X_{\pi(j)} - \bar{X}_n), \qquad t \in [0, 1].$$

Analogues of Theorems 2.1 and 2.2 can easily be established by exploiting the joint convergence of the partial sum processes $\sum_{j=1}^{\lfloor nt \rfloor} X_j$ and $\sum_{j=1}^{\lfloor nt \rfloor} |X_j|^\nu$. Theorems 2.3 and 2.4 also remain true (with some modifications in the corresponding limit processes) so that, conditionally on $\mathbf{X}$, the permuted sequence $\{A_{n,\pi}^{(\nu)}(t)\}$ cannot converge weakly in $\mathcal{D}[0, 1]$.

# 3. Proofs

**Proof of Theorem 2.2.** Note that, by the assumptions on $\{X_j\}$,

$$\{A_n(t) : t \in [0,1]\} \stackrel{\mathcal{D}}{=} \{A_{n,\pi}(t) : t \in [0,1]\},$$

so the result follows from Corollary 2.1.  □

From now on, we fix a realization of $\mathbf{X}$. Permuting the $X$'s is equivalent to selecting $n$ elements from $X_1, X_2, \ldots, X_n$ without replacement. The proof of Theorem 2.3 will be based on the order statistics $X_{1,n} \leq \cdots \leq X_{n,n}$ of $X_1, \ldots, X_n$. Our first goal is therefore to express $A_{n,\pi}(t)$ in terms of $X_{1,n}, \ldots, X_{n,n}$. To do so, let

$$\varepsilon_j^{(n)}(t) = \begin{cases} 1, & \text{if } X_{j,n} \text{ is among the first } \lfloor nt \rfloor \text{ elements chosen,} \\ 0, & \text{otherwise,} \end{cases}$$

that is, the sequence $\{\varepsilon_j^{(n)}(t)\}$ identifies those of the order statistics which are selected by the permutation $\pi$ in the first $\lfloor nt \rfloor$ positions. If there are identical observations, then we add multiples of $1/n^{2+1/\alpha}$ to the $X_{j,n}$'s to break the ties. Since this procedure will not change the asymptotics, from now on, we assume, without loss of generality, that all observations are different. It is then easy to see that

$$\sum_{j=1}^{\lfloor nt \rfloor} (X_{\pi(j)} - \bar{X}_n) = \sum_{j=1}^{n} (X_{j,n} - \bar{X}_n)\varepsilon_j^{(n)}(t) = \sum_{j=1}^{n} (X_{j,n} - \bar{X}_n)\bar{\varepsilon}_j^{(n)}(t),$$

where $\bar{\varepsilon}_j^{(n)}(t) = \varepsilon_j^{(n)}(t) - \mathrm{E}_{\mathbf{X}}\varepsilon_j^{(n)}(t) = \varepsilon_j^{(n)}(t) - \lfloor nt \rfloor / n$ is the centered version of $\varepsilon_j^{(n)}(t)$.

The proof of Theorem 2.3 will be the consequence of a series of lemmas to be stated next. Since $\{X_j\}$ is a sequence of random variables in the domain of attraction of a strictly $\alpha$-stable random variable with $\alpha \in (0,2)$, only the very small and very large order statistics will contribute asymptotically to $A_{n,\pi}(t)$, which is shown in Lemmas 3.1 and 3.3.

Let $\mathrm{Var}_{\mathbf{X}}$ denote conditional variance with respect to $\mathbf{X}$.

**Lemma 3.1.** *If (1.3) holds, then*

$$\lim_{K \to \infty} \limsup_{n \to \infty} P\left\{ P_{\mathbf{X}}\left\{ \frac{1}{n^{1/\alpha}L(n)} \left| \sum_{j=K+1}^{n-K} (X_{j,n} - \bar{X}_n)\varepsilon_j^{(n)}(t) \right| \geq \varepsilon \right\} \geq \delta \right\} = 0$$

*for all $\varepsilon > 0$ and $\delta > 0$.*

**Proof.** Observe that the statement of our theorems does not change if we replace $X_j$ by $X_j - c$ and thus we can assume that $\mu_\alpha = 0$ for all $\alpha \in (0, 2)$. Clearly $\mathrm{E}_{\mathbf{X}} \bar{\varepsilon}_j^{(n)}(t) = 0$ and

$$\mathrm{E}_{\mathbf{X}}[(\bar{\varepsilon}_j^{(n)}(t))^2] = \frac{\lfloor nt \rfloor}{n} - \left(\frac{\lfloor nt \rfloor}{n}\right)^2, \qquad j = 1, \ldots, n, \tag{3.1}$$

and, for $j \neq k$,

$$\mathrm{E}_{\mathbf{X}}[\bar{\varepsilon}_j^{(n)}(t) \bar{\varepsilon}_k^{(n)}(t)] = \frac{\lfloor nt \rfloor(\lfloor nt \rfloor - 1)}{n(n-1)} - \left(\frac{\lfloor nt \rfloor}{n}\right)^2 = -\frac{\lfloor nt \rfloor(n - \lfloor nt \rfloor)}{n^2(n-1)}. \tag{3.2}$$

Further we have

$$\mathrm{E}_{\mathbf{X}}\left[\sum_{j=K+1}^{n-K} (X_{j,n} - \bar{X}_n)\bar{\varepsilon}_j^{(n)}(t)\right] = 0$$

and

$$\mathrm{Var}_{\mathbf{X}}\left(\sum_{j=K+1}^{n-K} (X_{j,n} - \bar{X}_n)\bar{\varepsilon}_j^{(n)}(t)\right)$$

$$= \sum_{j,k=K+1}^{n-K} (X_{j,n} - \bar{X}_n)(X_{k,n} - \bar{X}_n)\mathrm{E}[\bar{\varepsilon}_j^{(n)}(t)\bar{\varepsilon}_k^{(n)}(t)]$$

$$\leq \frac{\lfloor nt \rfloor}{n} \frac{n - \lfloor nt \rfloor}{n} \sum_{j=K+1}^{n-K} (X_{j,n} - \bar{X}_n)^2 + \frac{1}{n-1}\left|\sum_{j,k=K+1, k\neq j}^{n-K} (X_{j,n} - \bar{X}_n)(X_{k,n} - \bar{X}_n)\right|$$

$$\leq \sum_{j=K+1}^{n-K} (X_{j,n} - \bar{X}_n)^2 + \frac{1}{n-1}\left|\left(\sum_{j=K+1}^{n-K} (X_{j,n} - \bar{X}_n)\right)^2 - \sum_{j=K+1}^{n-K} (X_{j,n} - \bar{X}_n)^2\right|$$

$$\leq 2\sum_{j=K+1}^{n-K} (X_{j,n} - \bar{X}_n)^2 + \frac{1}{n-1}\left(\sum_{j=K+1}^{n-K} (X_{j,n} - \bar{X}_n)\right)^2.$$

Now,

$$\sum_{j=K+1}^{n-K} (X_{j,n} - \bar{X}_n)^2 \leq \sum_{j=K+1}^{n-K} X_{j,n}^2 + 2|\bar{X}_n|\left|\sum_{j=K+1}^{n-K} X_{j,n}\right| + n\bar{X}_n^2.$$

Using (1.3) with $\mu_\alpha = 0$, we obtain

$$\frac{n\bar{X}_n^2}{n^{2/\alpha}L^2(n)} = \frac{1}{n}\left(\frac{1}{n^{1/\alpha}L(n)}\sum_{j=1}^{n} X_j\right)^2 = o_P(1) \qquad (n \to \infty).$$

As a consequence of Theorem 4.1 in Csörgő *et al.* (1986), we obtain that, for all $K \geq 1$,

$$\frac{1}{n^{1/\alpha}L(n)}\left|\sum_{j=K+1}^{n-K}X_{j,n}\right| = \mathcal{O}_P(1)$$

as $n \to \infty$. Therefore,

$$\frac{1}{n^{2/\alpha}L^2(n)}|\bar{X}_n|\left|\sum_{j=K+1}^{n-K}X_{j,n}\right| = \frac{1}{n}\left|\frac{1}{n^{1/\alpha}L(n)}\sum_{j=1}^{n}X_j\right|\left|\frac{1}{n^{1/\alpha}L(n)}\sum_{j=K+1}^{n-K}X_{j,n}\right| = o_P(1).$$

Clearly, the random variables $X_j^2$ belong to the domain of attraction of a stable law with parameter $\alpha/2$, with norming factor $n^{2/\alpha}L^2(n)$ and centering factor 0. Hence, Corollary 3.1, Theorem 4.1 and Proposition A.3 of Csörgő *et al.* (1986) imply that

$$\frac{1}{n^{2/\alpha}L^2(n)}\sum_{j=K+1}^{n-K}X_{j,n}^2 \xrightarrow{\mathcal{D}} A(K) \qquad (n \to \infty),$$

where

$$A(K) = q^{2/\alpha}\sum_{j=K+1}^{\infty}\frac{1}{S_j^{2/\alpha}} + p^{2/\alpha}\sum_{j=K+1}^{\infty}\frac{1}{(S_j^*)^{2/\alpha}}.$$

Clearly, $A(K) \to 0$ a.s. as $K \to \infty$. Another application of Theorem 4.1 in Csörgő *et al.* (1986) yields that, for any $K \geq 1$ and $n \to \infty$,

$$\frac{1}{n-1}\frac{1}{n^{2/\alpha}L^2(n)}\left(\sum_{j=K+1}^{n-K}(X_{j,n}-\bar{X}_n)\right)^2 = \frac{1}{n-1}\left(\frac{1}{n^{1/\alpha}L(n)}\sum_{j=K+1}^{n-K}(X_{j,n}-\bar{X}_n)\right)^2 = o_P(1).$$

The assertion of Lemma 3.1 is now obtained from Markov's inequality.                    □

Without loss of generality, we henceforth assume that all processes used in this paper are defined on the same probability space. In the following, we shall utilize a result of Berkes and Philipp (1979) to show that, on the remaining index range $j \in [0, K] \cup [n-K+1, n]$, the dependent random variables $\{\varepsilon_j^{(n)}(t)\}$ can be approximated with the sequence of independent Bernoulli variables $\{\delta_j^{(n)}(t)\}$ introduced above.

**Lemma 3.2.** *If (1.3) holds, then, for each $n$ and each $t \in (0,1)$, there exist independent, identically distributed random variables $\delta_j^{(n)}(t)$, $j = 1, \ldots, n$, independent of $\{X_j\}$, with*

$$P\{\delta_j^{(n)}(t) = 1\} = \frac{\lfloor nt \rfloor}{n} \quad and \quad P\{\delta_j^{(n)}(t) = 0\} = \frac{n - \lfloor nt \rfloor}{n} \tag{3.3}$$

*such that*

$$P_{\mathbf{X}}\left\{\sum_{j=1}^{K}(X_{j,n}-\bar{X}_n)(\varepsilon_j^{(n)}(t)-\delta_j^{(n)}(t))\neq 0\right\}\leq\frac{48K^2}{n} \qquad (3.4)$$

*and*

$$P_{\mathbf{X}}\left\{\sum_{j=n-K+1}^{n}(X_{j,n}-\bar{X}_n)(\varepsilon_j^{(n)}(t)-\delta_j^{(n)}(t))\neq 0\right\}\leq\frac{48K^2}{n} \qquad (3.5)$$

*for all $K=1,\ldots,\lfloor n/2\rfloor$.*

Note that relations (3.4) and (3.5) are obvious for $K\geq\sqrt{n/48}$, but we will use the lemma for constant $K$.

**Proof of Lemma 3.2.** For $j\geq 1$, let $\gamma_{2j-1}=\varepsilon_j^{(n)}(t)$ and $\gamma_{2j}=\varepsilon_{n-j+1}^{(n)}(t)$. To replace the dependent $\varepsilon_j^{(n)}(t)$ with the independent $\delta_j^{(n)}(t)$ given in (3.3), we need to derive an upper bound for the difference between the conditional probability $P\{\gamma_k=a_k|\gamma_{k-1}=a_{k-1},\ldots,\gamma_1=a_1\}$ and the probability $P\{\gamma_k=a_k\}$. To do so, consider a set $A=\{x_1,\ldots,x_n\}$ and choose from it a subset of $\lfloor nt\rfloor$ elements. The probability that this subset contains $r$ fixed elements of $A$ but does not contain another fixed $s$ elements of $A$ is $p_{r,s}^{(n)}(t)=\binom{n-r-s}{\lfloor nt\rfloor-r}/\binom{n}{\lfloor nt\rfloor}$. Provided that $k=r+s$, it holds that $P\{\gamma_1=a_1,\ldots,\gamma_k=a_k\}=p_{r,s}^{(n)}(t)$, where $r$ of the coefficients $a_\nu$ are equal to 1 and $s$ are equal to 0. Similarly, $P\{\gamma_1=a_1,\ldots,\gamma_k=a_k,\gamma_{k+1}=1\}=p_{r+1,s}^{(n)}(t)$. Consequently,

$$P\{\gamma_{k+1}=1|\gamma_1=a_1,\ldots,\gamma_k=a_k\}=\frac{p_{r+1,s}^{(n)}(t)}{p_{r,s}^{(n)}(t)}=\frac{\lfloor nt\rfloor-r}{n-r-s}.$$

Since it can be assumed that $1\leq k\leq n/2$, we can estimate

$$\left|\frac{\lfloor nt\rfloor-r}{n-r-s}-t\right|\leq\frac{1+r+s}{n-r-s}=\frac{1+k}{n(1-k/n)}<\frac{1+k}{n}\left(1+\frac{2k}{n}\right)\leq\frac{4k}{n}.$$

In a similar fashion, we can obtain

$$|P\{\gamma_{k+1}=1\}-t|\leq\frac{4k}{n}.$$

Thus, for $a_{k+1}=1$ (and consequently for $a_{k+1}=0$), we have

$$|P\{\gamma_{k+1}=a_{k+1}|\gamma_k=a_k,\ldots,\gamma_1=a_1\}-P\{\gamma_{k+1}=a_{k+1}\}|\leq\frac{8k}{n}$$

for all $k = 2, \ldots, n$. Hence, applying Theorem 2 of Berkes and Philipp (1979), there exist independent Bernoulli random variables $\delta_j^{(n)}(t)$, $j = 1, \ldots, n$, satisfying (3.3) such that

$$P\left\{ |\varepsilon_k^{(n)}(t) - \delta_k^{(n)}(t)| \geq \frac{48k}{n} \right\} \leq \frac{48k}{n}$$

and

$$P\left\{ |\varepsilon_{n-k+1}^{(n)}(t) - \delta_{n-k+1}^{(n)}(t)| \geq \frac{48k}{n} \right\} \leq \frac{48k}{n}$$

for $k = 1, \ldots, \lfloor n/2 \rfloor$. Since the variables $\varepsilon_j^{(n)}(t)$ and $\delta_j^{(n)}(t)$ take only the values 0 and 1, the last two formulas imply that, with the exception of a set of probability not exceeding $48K^2/n$, all differences $\varepsilon_j^{(n)}(t) - \delta_j^{(n)}(t)$ in the sum in (3.4), and similarly in the sum in (3.5), are equal to 0. □

The following lemma is the counterpart of Lemma 3.1, in which the $\bar{\varepsilon}_j^{(n)}(t)$ are replaced by mean-corrected variables $\bar{\delta}_j^{(n)}(t)$.

**Lemma 3.3.** *If (1.3) holds, then*

$$\lim_{K \to \infty} \limsup_{n \to \infty} P\left\{ P_{\mathbf{X}}\left\{ \frac{1}{n^{1/\alpha} L(n)} \left| \sum_{j=K+1}^{n-K} (X_{j,n} - \bar{X}_n) \bar{\delta}_j^{(n)}(t) \right| \geq \varepsilon \right\} \geq \delta \right\} = 0$$

*for all $\varepsilon > 0$ and $\delta > 0$, where $\bar{\delta}_j^{(n)}(t) = \delta_j^{(n)}(t) - \lfloor nt \rfloor / n$.*

**Proof.** The assertion follows along the lines of the proof of Lemma 3.1. □

Let $S_{j,n} = (X_{j,n} - \bar{X}_n)/T_n$ denote the standardized extreme order statistics. The joint weak convergence of the $S_{j,n}$, $j \in [1, K] \cup [n - K + 1, n]$, is given in the following lemma. Recall from Section 2 that $\{S_j : j \geq 1\}$ and $\{S_j^* : j \geq 1\}$ are sums of independent sequences consisting of independent exponential random variables with parameter 1.

**Lemma 3.4.** *If (1.3) holds, then, as $n \to \infty$,*

$$(S_{j,n} : j \in [1, K] \cup [n - K + 1, n]) \overset{\mathcal{D}}{\longrightarrow} \frac{1}{M}\left( -\frac{w_1}{S_1^{1/\alpha}}, \ldots, -\frac{w_1}{S_K^{1/\alpha}}, \frac{w_2}{(S_K^*)^{1/\alpha}}, \ldots, \frac{w_2}{(S_1^*)^{1/\alpha}} \right),$$

*where $w_1 = q^{1/\alpha}$, $w_2 = p^{1/\alpha}$ and $M$ is the random variable defined in (2.5).*

**Proof.** Csörgő *et al.* (1986), page 110, showed that

$$\frac{1}{n^{1/\alpha} L(n)} (X_{j,n} : j \in [1, K] \cup [n - K + 1, n])$$

$$\xrightarrow{\mathcal{D}} \left( -\frac{w_1}{S_1^{1/\alpha}}, \dots, -\frac{w_1}{S_K^{1/\alpha}}, \frac{w_2}{(S_K^*)^{1/\alpha}}, \dots, \frac{w_2}{(S_1^*)^{1/\alpha}} \right).$$

Since, by (1.3), $\bar{X}_n/T_n \xrightarrow{P} 0$ and, clearly, $T_n = \max\{|X_{1,n}|, |X_{n,n}|\}$, Lemma 3.4 is proved. $\qquad\square$

**Proof of Theorem 2.3.** In view of Lemma 3.1, we have

$$\limsup_{n\to\infty} \mathcal{L}\left( \mathrm{dist}_{\mathbf{X}} \frac{1}{n^{1/\alpha}L(n)} \sum_{j=1}^n (X_{j,n} - \bar{X}_n) \bar{\varepsilon}_j^{(n)}(t), \right.$$

$$\left. \mathrm{dist}_{\mathbf{X}} \frac{1}{n^{1/\alpha}L(n)} \sum_{j \in [1,K] \cup [n-K+1,n]} (X_{j,n} - \bar{X}_n) \bar{\varepsilon}_j^{(n)}(t) \right) = B(K),$$

where $B(K) \xrightarrow{P} 0$ as $K \to \infty$. Here, $\mathrm{dist}_{\mathbf{X}}$ denotes conditional distribution with respect to $\mathbf{X}$ and $\mathcal{L}$ is the Lévy distance. We used here the fact that if, for two random variables $\xi$ and $\eta$, we have $P(|\xi - \eta| \geq \varepsilon) \leq \varepsilon$, then the Lévy distance of the distributions of $\xi$ and $\eta$ is not greater than $\varepsilon$. By Lemma 3.3, the previous relation remains valid if we replace $\bar{\varepsilon}_j^{(n)}(t)$ by $\bar{\delta}_j^{(n)}(t)$. Further, $\bar{\varepsilon}_j^{(n)}(t) - \bar{\delta}_j^{(n)}(t) = \varepsilon_j^{(n)}(t) - \delta_j^{(n)}(t)$ and thus Lemma 3.2 implies that for any fixed $K$, the Lévy distance of the conditional distributions

$$\mathrm{dist}_{\mathbf{X}} \frac{1}{n^{1/\alpha}L(n)} \sum_{j \in [1,K] \cup [n-K+1,n]} (X_{j,n} - \bar{X}_n) \bar{\varepsilon}_j^{(n)}(t)$$

and

$$\mathrm{dist}_{\mathbf{X}} \frac{1}{n^{1/\alpha}L(n)} \sum_{j \in [1,K] \cup [n-K+1,n]} (X_{j,n} - \bar{X}_n) \bar{\delta}_j^{(n)}(t)$$

tends to 0 in probability as $n \to \infty$. Since $T_n/(n^{1/\alpha}L(n))$ has a non-degenerate limit distribution, the above statements remain valid if we replace the norming factor $n^{1/\alpha}L(n)$ by $T_n$. Thus, it suffices to consider the limiting behavior of

$$\mathrm{dist}_{\mathbf{X}} \frac{1}{T_n} \left( \sum_{j \in [1,K] \cup [n-K+1,n]} (X_{j,n} - \bar{X}_n) \bar{\delta}_j^{(n)}(t) \right) = \mathrm{dist}_{\mathbf{X}} \left( \sum_{j \in [1,K] \cup [n-K+1,n]} S_{j,n} \bar{\delta}_j^{(n)}(t) \right)$$

for fixed $K$. By the Skorokhod–Dudley–Wichura representation theorem (see, e.g., Shorack and Wellner (1986), page 47), on a sufficiently large probability space, one can redefine the vectors $(S_{j,n} : j \in [1, K] \cup [n - K + 1, n])$ and the independent sequences $\{S_j : j \geq 1\}$, $\{S_j^* : j \geq 1\}$ such that the convergence relation in Lemma 3.4 holds almost surely. Since the Lévy distance of the distribution of two sums

$$\sum_{j \in [1,K] \cup [n-K+1,n]} c_j \bar{\delta}_j^{(n)}(t) \quad \text{and} \quad \sum_{j \in [1,K] \cup [n-K+1,n]} c_j' \bar{\delta}_j^{(n)}(t)$$

is not greater than $\sum_{j \in [1,K] \cup [n-K+1,n]} |c_j - c'_j|$ for any real sequences $\{c_j\}, \{c'_j\}$, it follows that for the redefined variables $S_{j,n}, S_j, S_j^*$, we have

$$\mathcal{L}\Bigg(\mathrm{dist}_{\mathbf{S}} \sum_{j \in [1,K] \cup [n-K+1,n]} S_{j,n} \bar{\delta}_j^{(n)}(t),$$

$$\mathrm{dist}_{\mathbf{S}} \frac{1}{M} \left[ -w_1 \sum_{j=1}^{K} \frac{1}{S_j^{1/\alpha}} \bar{\delta}_j^{(n)}(t) + w_2 \sum_{j=1}^{K} \frac{1}{(S_j^*)^{1/\alpha}} \bar{\delta}_{n-j+1}^{(n)}(t) \right] \Bigg) \xrightarrow{P} 0$$

as $n \to \infty$. Clearly, this redefinition will change the random distribution

$$\mathrm{dist}_{\mathbf{X}} \sum_{j \in [1,K] \cup [n-K+1,n]} S_{j,n} \bar{\delta}_j^{(n)}(t)$$

(as it will be defined on a new probability space), but not its distribution. Hence, it suffices to show that the Lévy distance of the conditional distributions

$$\mathrm{dist}_{\mathbf{S}} \frac{1}{M} \left( -w_1 \sum_{j=1}^{K} \frac{1}{S_j^{1/\alpha}} \bar{\delta}_j^{(n)}(t) + w_2 \sum_{j=1}^{K} \frac{1}{(S_j^*)^{1/\alpha}} \bar{\delta}_{n-j+1}^{(n)}(t) \right) \qquad (3.6)$$

and

$$\mathrm{dist}_{\mathbf{S}} \frac{1}{M} \left( -w_1 \sum_{j=1}^{K} \frac{1}{S_j^{1/\alpha}} (\delta_j(t) - t) + w_2 \sum_{j=1}^{K} \frac{1}{(S_j^*)^{1/\alpha}} (\delta_j^*(t) - t) \right) \qquad (3.7)$$

tends to 0 a.s. for any fixed $K \geq 1$ and that the conditional distribution in (3.7) converges a.s. to the same expression with $K = \infty$. The first statement is obvious from the fact that $\bar{\delta}_j^{(n)}(t) = \delta_j(t) - \lfloor nt \rfloor / n$ and the second statement follows from the fact that, by the strong law of large numbers, we have $S_j/j \to 1$ a.s. and thus

$$\mathrm{Var}_{\mathbf{S}} \left( \frac{1}{S_j^{1/\alpha}} (\delta_j(t) - t) \right) = \mathcal{O}(j^{-2/\alpha}) \qquad \text{a.s.},$$

whence it follows that the series

$$\sum_{j=1}^{\infty} \frac{1}{S_j^{1/\alpha}} (\delta_j(t) - t) \quad \text{and} \quad \sum_{j=1}^{\infty} \frac{1}{(S_j^*)^{1/\alpha}} (\delta_j^*(t) - t)$$

are, conditionally on $\mathbf{S}$, a.s. convergent by the Kolmogorov two series theorem. Note, finally, that, by a theorem of Lévy (1931), page 150 (see also Breiman (1968), page 51, Problem 16), the distribution

$$\mathrm{dist}_{\mathbf{S}} \frac{1}{M} \left( -w_1 \sum_{j=1}^{\infty} \frac{1}{S_j^{1/\alpha}} (\delta_j(t) - t) + w_2 \sum_{j=1}^{\infty} \frac{1}{(S_j^*)^{1/\alpha}} (\delta_j^*(t) - t) \right)$$

is a.s. continuous. Since weak convergence of (ordinary) distributions to a continuous limit implies the pointwise convergence of the corresponding distribution functions, Theorem 2.3 follows. □

The proof of Theorem 2.4 requires a simple modification of Lemma 3.2.

**Lemma 3.5.** *If (1.3) holds, then, for each $n$ and each choice of $0 < t_1 < \cdots < t_r < 1$, there exist independent, identically distributed random vectors $(\delta_j^{(n)}(t_1), \ldots, \delta_j^{(n)}(t_r))$, $j = 1, \ldots, n$, independent of $\{X_j\}$, such that*

$$(\delta_j^{(n)}(t_1), \ldots, \delta_j^{(n)}(t_r)) \stackrel{\mathcal{D}}{=} \left( I\left\{ U \le \frac{\lfloor nt_1 \rfloor}{n} \right\}, \ldots, I\left\{ U \le \frac{\lfloor nt_r \rfloor}{n} \right\} \right),$$

*where $U$ denotes a uniform random variable on $[0,1]$ and, further, for all $\delta > 0$,*

$$P_{\mathbf{X}}\left\{ \max_{1 \le \ell \le r} \frac{1}{n^{1/\alpha} L(n)} \left| \sum_{j=1}^{K} (X_{j,n} - \bar{X}_n)(\varepsilon_j^{(n)}(t_\ell) - \delta_j^{(n)}(t_\ell)) \right| > \delta \right\} \stackrel{P}{\longrightarrow} 0$$

*and*

$$P_{\mathbf{X}}\left\{ \max_{1 \le \ell \le r} \frac{1}{n^{1/\alpha} L(n)} \left| \sum_{j=n-K+1}^{n} (X_{j,n} - \bar{X}_n)(\varepsilon_j^{(n)}(t_\ell) - \delta_j^{(n)}(t_\ell)) \right| > \delta \right\} \stackrel{P}{\longrightarrow} 0.$$

**Proof.** The assertion of Lemma 3.5 follows by applying the approximation procedure in the proof of Lemma 3.2 to the vector $(\varepsilon_j^{(n)}(t_1), \ldots, \varepsilon_j^{(n)}(t_r))$ instead of $\varepsilon_j^{(n)}(t)$. Since the changes are routine, we omit the details. □

**Proof of Theorem 2.4.** The proof can be given by repeating the arguments developed in the proof of Theorem 2.3, replacing Lemma 3.2 with Lemma 3.5. □

# Acknowledgements

Research partially supported by NATO grant PST.EAP.CLG 980599, NSF Grant DMS-06-4670, Grant RGC-HKUST 6428/06H and Hungarian National Foundation for Scientific Research Grants T 043037 and K 61052.